# Constrained Spacecraft Relative Motion Planning Exploiting Periodic Natural Motion Trajectories and Invariance

Gregory R. Frey[*]
*University of Michigan, Ann Arbor, Michigan 48109-2140*
Christopher D. Petersen[†]
*U.S. Air Force Research Laboratory, Kirtland Air Force Base, Albuquerque, New Mexico, 87117-5776*
Frederick A. Leve[‡]
*Air Force Office of Scientific Research, Arlington, Virginia 22203*
Ilya V. Kolmanovsky[§] and Anouck R. Girard[¶]
*University of Michigan, Ann Arbor, Michigan 48109-2140*

## Abstract

**Spacecraft relative motion planning is concerned with the design and execution of maneuvers relative to a nominal target. These types of maneuvers are frequently utilized in missions such as rendezvous and docking, satellite inspection and formation flight where exclusion zones representing spacecraft or other obstacles must be avoided. The presence of these exclusion zones leads to non-linear and non-convex constraints which must be satisfied. In this paper, a novel approach to spacecraft relative motion planning with obstacle avoidance and thrust constraints is developed. This approach is based on a graph search applied to a virtual net of closed (periodic) natural motion trajectories, where the natural motion trajectories represent virtual net nodes (vertices), and adjacency and connection information is determined by conditions defined in terms of safe, positively-invariant tubes built around each trajectory. These conditions guarantee that transitions from one natural motion trajectory to another natural motion trajectory can be completed without constraint violations. The proposed approach improves the flexibility of a previous approach based on the use of forced equilibria, and has other advantages in terms of reduced fuel consumption and passive safety. The resulting maneuvers, if planned on-board, can be executed directly or, if planned off-board, can be used to warm start trajectory optimizers to generate further improvements.**

# Nomenclature

| | | |
|---|---|---|
| $A$, $A_c$, $\bar{A}$ | = | Discrete-time, continuous-time and closed-loop dynamics matrices |
| $B$ | = | Discrete-time input matrix |
| $e$ | = | State error |
| $J$ | = | Trajectory cost |
| $k$ | = | Discrete-time instant (integer) |
| $K$ | = | State-feedback gain matrix |
| $P$ | = | Positive-definite ellipsoidal shape matrix |
| $u$ | = | Control vector |
| $u_{max}$ | = | Maximum allowable control |
| $X$ | = | Spacecraft state vector consisting of relative positions and velocities, $x$, $y$, $z$, $\dot{x}$, $\dot{y}$, $\dot{z}$ |
| $X_n$ | = | State vector along a natural motion trajectory |
| $X_{ni}$ | = | State vector along natural motion trajectory $\mathcal{N}_i$ |
| $\delta$ | = | Integer corresponding to initial controller reference point along a natural motion trajectory |
| $\Delta T$ | = | Discrete-time update period |

[*]Ph.D. Candidate, Department of Aerospace Engineering, gfrey@umich.edu, AIAA Member
[†]Research Engineer, Guidance, Navigation and Control Group, christopher.petersen.16@us.af.mil, AIAA Member
[‡]Program Officer, Dynamics and Controls Program, frederick.leve@us.af.mil, AIAA Member
[§]Professor, Department of Aerospace Engineering, ilya@umich.edu, AIAA Member
[¶]Associate Professor, Department of Aerospace Engineering, anouck@umich.edu, AIAA Member

| | | |
|---|---|---|
| $\Xi, \Xi_w$ | = | Unweighted and weighted connection arrays |
| $\Pi, \Pi_w$ | = | Unweighted and weighted adjacency matrices |
| $\rho_k^s$ | = | Ellipsoidal scale factors for safe sets |
| $\rho_k$ | = | Ellipsoidal scale factors used to generate safe, positively-invariant tubes |
| $\rho_u, \rho_{Oi,k}$ | = | Maximum possible ellipsoidal scale factors considering control constraints, or the $i$th exclusion zone constraint |
| $\mathcal{B}(Z, \gamma)$ | = | Ball of radius $\gamma$ centered at state vector $Z$ |
| $\mathcal{E}_{k,\mathcal{N}}$ | = | Ellipsoidal set centered at $X_n(k)$ along natural motion trajectory $\mathcal{N}$ with scale factor $\rho_k$ |
| $\mathcal{E}_{k,\mathcal{N}}^s$ | = | Safe ellipsoidal set centered at $X_n(k)$ along natural motion trajectory $\mathcal{N}$ with scale factor $\rho_k^s$ |
| $\mathcal{N}$ | = | Set of state vectors corresponding to a closed natural motion trajectory |
| $\mathcal{O}(s_i, S_i)$ | = | Ellipsoidal exclusion zone centered at point $s_i$ with shape matrix $S_i$ |
| $\mathcal{T}_\mathcal{N}^s$ | = | Safe tube for natural motion trajectory $\mathcal{N}$ |
| $\mathcal{T}_\mathcal{N}$ | = | Safe, positively invariant tube for natural motion trajectory $\mathcal{N}$ |
| $\mathbb{R}$ | = | Set of real numbers |
| $\mathbb{Z}$ | = | Set of integers |

# 1 Introduction

Relative motion planning must frequently account for obstacles, represented by exclusion zones, in order to ensure safe operations. For a satellite mission, these obstacles may be pieces of orbital debris, other spacecraft, or areas which must be avoided due to sensor constraints. Obstacle avoidance requirements often result in non-linear and/or non-convex constraints on vehicle motion, complicating the application of conventional trajectory optimization methods. While the problem of motion planning with obstacle avoidance is commonly encountered in many fields such as robotics [1], several factors make the application to spacecraft motion planning unique. Firstly, fuel efficiency is important for spacecraft as refueling is not possible. Secondly, spacecraft frequently have limited on-board computing capabilities, thereby requiring fast and efficient on-board trajectory computation algorithms. Finally, spacecraft dynamics include periodic behavior and many Natural Motion Trajectories (NMTs) which can be followed with no control usage and utilized to generate fuel efficient trajectories.

An extensive body of literature exists related to the problem of spacecraft motion planning with obstacle avoidance. A concise review of this literature is provided here to highlight the range of methods that have been applied to the problem, and to motivate the proposed approach. A method for calculating fuel-optimal trajectories with obstacle avoidance is formulated in [2] using mixed-integer linear programming. A method to generate delta-V ($\Delta V$) optimal paths to inspect another spacecraft while avoiding keep out zones is described in [3], using Sparse Optimal Control Software (SPOCS) [4]. The generation of passively safe paths, i.e., trajectories that guarantee collision avoidance in the presence of anomalous behaviors such as thruster failure, is considered in [5] using receding horizon control.

Spacecraft formation reconfiguration while avoiding collisions has been proposed using several methods including state-constrained optimal control techniques [6], a passivity-based sliding surface controller [7], and heuristics involving separation planes [8]. An eccentricity/inclination vector separation method to ensure adequate separation distances between spacecraft in formation is described in [9]. This method has been studied for use in both the GRACE and PRISMA formation flight missions [10], [11]. Methods for trajectory planning with obstacle avoidance using artificial potential functions have also been considered in, e.g., [12], [13] and [14]. Approaches based on solving non-convex trajectory optimization problems with a sequence of convex optimization problems have been proposed in [15] and [16].

Motion planning using graph search is desirable in spacecraft applications because the efficiency and simplicity of certain algorithms, such as Dijkstra's algorithm [17], make implementation on-board a satellite with limited computational capability achievable. A graph theoretic framework is applied in [18] to Leader Following (LF) spacecraft formation control. Rapidly Exploring Random Tree (RRT) and similar algorithms have also been applied to spacecraft motion planning while accounting for exclusion zones [19]. A Fast Marching Tree (FMT) algorithm was applied in [20] to develop safe paths for satellite rendezvous.

A framework for spacecraft relative motion trajectory planning with obstacle avoidance that exploits graph search on a virtual net consisting of static points (forced equilibria) in Hill's relative motion frame, [21], has been proposed in [22] (see also [23] for more recent work). In [22], safe (constraint admissible) positively-invariant sets were used to determine feasibility of node-to-node transitions. Furthermore, it has been shown that the approach can be easily extended to include bounded disturbances and moving obstacles. In this paper, we demonstrate that it is possible to integrate closed non-equilibrium NMTs into this framework. Specifically, NMTs are used to represent virtual net nodes, and adjacency is determined by conditions defined in terms of safe, positively invariant tubes around each trajectory. The utilization of closed NMTs has several advantages compared to the forced equilibria considered in [22]. Firstly, travelling along closed NMTs in steady-state is possible with zero fuel consumption, while zero fuel consumption is only achieved for forced equilibria along the in-track axis. Secondly, the use of closed NMTs expands the set of trajectories available to compose the overall maneuver from, while ensuring the resulting maneuvers are fuel efficient, i.e., when the maneuver consists of NMTs and transfers connecting them, fuel is consumed only during the transfers and to compensate for perturbations. Thirdly, the use of closed NMTs has advantages in terms of passive safety as the spacecraft can remain on a closed NMT, which does not intersect known obstacles, and avoid collisions even if thrust is temporarily lost. Since closed NMTs are open-loop unstable, (however, not exponentially unstable), generally these passive safety properties can be exploited over short periods of time, after which thrust-based control must be regained.

For background information on invariance, safe positively-invariant sets, and their use, see, e.g., [24–27]. Invariant tubes are utilized in [28–30] to account for the effect of unmeasured disturbances. Other related work on trajectory planning with obstacle avoidance, not specifically developed for spacecraft, includes [31], which developed an LQR-Trees algorithm to exploit a set of trajectories, calculated using trajectory optimization algorithms, and stabilized using time-varying LQR controllers. The regions of attraction for these trajectories "probabilistically cover" the controllable state space. More recent related developments include [32], where a trajectory planning method is developed using invariant "funnels" around a set of open loop maneuvers. These funnels are used to piece together multiple trajectories, forming a path which avoids obstacles. The trajectory planning is accomplished on-line and can be re-computed during execution if additional obstacles are discovered. Finally, reference [33] develops a control law which can be applied to track a nominal trajectory, and uses this control law to form an invariant tube around the trajectory. The use of this control law guarantees that motion will remain within the tube, and thus constraints are satisfied. Because the controller can be applied to any trajectory, the nominal path may be adjusted during execution while keeping the control law unchanged. Our work is different from [31–33] in that it is focused on taking into account spacecraft relative motion dynamics, which are open-loop unstable, and the trajectories considered here are periodic NMTs which can be obtained without resorting to trajectory optimization methods, and which can be followed with zero fuel consumption (or minimal fuel consumption if perturbations are considered) once reached. Additionally, fixed gain LQR controllers are used in our work which also leads to simple implementation. Finally, for trajectory planning, we use simple graph search on a virtual net with node adjacency rules that already account for known obstacles, hence feasible trajectories (satisfying both control and exclusion zone constraints) may be planned with minimal computations.

In our approach, the adjacency of nodes representing closed NMTs is determined by forming safe, positively-invariant tubes around each NMT. These tubes are generated as unions of safe ellipsoidal sets centered at points along each NMT. Within each tube, constraints on both control and state variables are satisfied, and thus constraint satisfaction, including obstacle avoidance, is guaranteed for any trajectory which stays within the tubes. In contrast to [22], where the positive invariance of ellipsoidal sets around forced equilibria was guaranteed regardless of the set size, in this work, the size of each ellipsoidal set in the tube around a closed NMT must be appropriately selected to ensure positive-invariance. Two methods of selecting the ellipsoidal set sizes are developed and proven to yield safe, positively-invariant tubes. One of these methods is conservative, i.e., forms a relatively small tube, but requires minimal computations. A second method forms the largest possible safe, positively-invariant tube consisting of ellipsoidal sets, at the expense of slightly increased computational load. This increase in tube size provides additional flexibility in trajectory planning.

Trajectories for the spacecraft to follow are generated by graph search using Dijkstra's algorithm, to produce a sequence of NMTs. The spacecraft traverses this sequence of NMTs using a fixed gain state-feedback control law with a time-varying reference along the current NMT. When the spacecraft reaches the prescribed transfer location, the controller reference is switched to the next NMT in the sequence. Trajectories with improved fuel efficiency are obtained by selecting appropriate costs for the node adjacency matrix. The on-board calculation of safe trajectories is facilitated by the introduction of a connection array which provides the starting point and initial controller reference point to be used to execute transfers between any two adjacent NMTs. These trajectories may be either executed as-is, or used to warm-start open-loop trajectory optimization algorithms.

As a summary, the specific contributions of this work are: a) development and utilization of a virtual net consisting of nodes corresponding to closed NMTs for planning spacecraft relative motion trajectories that can be closed-loop followed, b) formulation of two procedures to generate safe, positively invariant tubes around each NMT, c) generation of a connection array which can be used to simplify the on-line calculations needed to generate safe trajectories between pairs of NMTs, and d) demonstration of the ability of the proposed framework to generate trajectories which avoid obstacles through simulations.

This paper is organized as follows. In Section 2, the spacecraft model is summarized including the dynamics, control law, constraints and a description of the types of NMTs considered. Section 3 describes the generation of safe, positively-invariant tubes and provides two procedures which can be used to generate these tubes. Section 4 introduces the virtual net which is used to reduce the problem of trajectory planning to a conventional graph search. Methods are provided to determine both the adjacency of the virtual net and connection information which provides parameters used to generate safe transfers between NMTs. Simulation results are presented in Section 5 to illustrate these trajectory planning methods. Finally, Section 6 contains concluding remarks.

## 2 Preliminaries

### 2.1 Spacecraft Model

The spacecraft dynamics model is formulated in Hill's reference frame, which has the origin at a specified location on a nominal circular spacecraft orbit. The $x$-axis is in the radial direction, defined by the line from the center of the earth to the origin, the $z$-axis is in the direction of the nominal orbit angular momentum vector, and the $y$-axis completes the right-handed coordinate system. As the origin of this reference frame moves along the nominal circular orbit, the reference frame rotates at a rate equal to the mean motion of the circular orbit. Circular orbits are considered here as 79% of satellite orbits are nearly circular (have an eccentricity of less than 0.025 [34]), and this case yields time-invariant relative motion dynamics with closed NMTs that are easy to characterize.

In Hill's frame, the motion of a spacecraft relative to the origin is expressed using the linearized Clohessy-Wiltshire (CW) equations [35], which in discrete-time are

$$X(k+1) = AX(k) + Bu(k), \tag{1}$$

where $k \in \mathbb{Z}_{\geq 0}$ denotes the discrete-time instants,

$$X = [x \ y \ z \ \dot{x} \ \dot{y} \ \dot{z}]^T, \tag{2}$$

and where $x$, $y$, $z$, are the relative coordinates of the spacecraft in Hill's frame, $\dot{x}$, $\dot{y}$, and $\dot{z}$ are components of the relative velocity vector, and $u(k)$ is the control vector corresponding to continuous thrust forces.

Assuming an update period of $\Delta T$ seconds, the discrete-time dynamics and input matrices have the following form,

$$A = \exp(A_c \Delta T), \tag{3}$$

$$B = \int_0^{\Delta T} \exp(A_c(\Delta T - \tau))d\tau \begin{bmatrix} 0_{3\times 3} \\ \frac{1}{m}I_{3\times 3} \end{bmatrix}, \tag{4}$$

where $m$ is the mass of the spacecraft, $0_{3\times 3}$ is the $3 \times 3$ matrix consisting of all 0s, $I_{3\times 3}$ is the $3 \times 3$ identity

matrix, and

$$A_c = \begin{bmatrix} 0 & 0 & 0 & 1 & 0 & 0 \\ 0 & 0 & 0 & 0 & 1 & 0 \\ 0 & 0 & 0 & 0 & 0 & 1 \\ 3\omega^2 & 0 & 0 & 0 & 2\omega & 0 \\ 0 & 0 & 0 & -2\omega & 0 & 0 \\ 0 & 0 & -\omega^2 & 0 & 0 & 0 \end{bmatrix}, \quad (5)$$

where $\omega = \sqrt{\mu/R_0^3}$ is the mean motion of the nominal circular orbit, $\mu$ is Earth's gravitational parameter, and $R_0$ is the orbital radius of the nominal circular orbit.

## 2.2 Natural Motion Trajectories

A NMT is defined as a solution to (1) with $u = 0$. Depending on the initial condition, NMTs can take a variety of forms, including ellipses, spirals, lines and stationary points [36]. With an initial condition $X(0) = \bar{X}_0$ selected such that

$$\dot{y}(0) = -2\omega x(0), \quad (6)$$

the resulting trajectory will be periodic, i.e., closed, with a period equal to that of the nominal circular orbit, $\tau = 2\pi/\omega$ [37].

Closed NMTs can be stationary points along the $y$-axis (in-track), periodic line segments in the $y$-$z$ plane with $y(k) = y(0)$, or ellipses centered at a point along the $y$-axis. Figure 1 shows examples of these types of closed NMTs. Methods to generate initial conditions for these types of closed NMTs are available, see,

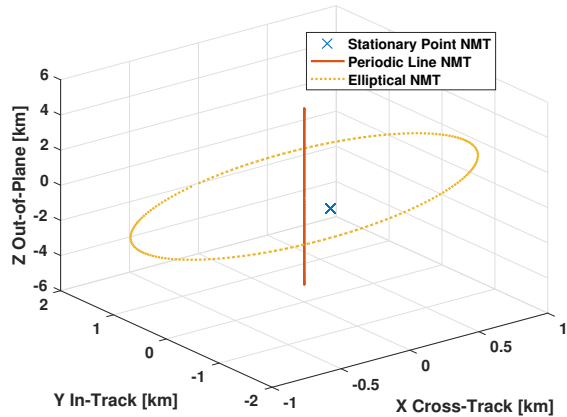

**Fig. 1** Examples of different types of closed NMTs

e.g., [38], [39], and are also included in Appendix A for completeness.

If $\Delta T$ is chosen such that $\frac{\tau}{\Delta T} \in \mathbb{Z}_{>0}$, where $\tau$ is the period of the nominal circular orbit, a closed NMT $\mathcal{N}$ starting from a specified initial condition $\bar{X}_0$ can be defined as a finite set of state vectors,

$$\mathcal{N}(\bar{X}_0) = \mathcal{N} = \left\{ X_n(k) \mid X_n(0) = \bar{X}_0, \ X_n(k+1) = AX_n(k), \ k \in [0, k_{max}] \right\}, \quad (7)$$

where $k_{max} = \frac{\tau}{\Delta T} - 1$.

**Remark 1:** Note that the set of state vectors defined by (7) completely defines the NMT because for $k > k_{max}$, the sequence of state vectors repeats, i.e., $X_n(k_{max} + 1) = X_n(0)$, $X_n(k_{max} + 2) = X_n(1)$, etc. In general, $X_n(k) = X_n(\tilde{k})$ where

$$\tilde{k} = \mathrm{mod}(k, \ k_{max} + 1), \quad (8)$$

and where the modulo function $\mathrm{mod}(x, y)$ returns the remainder after division of $x$ by $y$. In all subsequent developments, any index $k > k_{max}$ is taken to be the equivalent index $\tilde{k} \in [0, k_{max}]$ given by (8).

## 2.3 Spacecraft Control Law

The nominal feedback law which guides the spacecraft to a desired closed NMT is given by

$$u(k) = K(X(k) - X_n(k+\delta)), \tag{9}$$

where $K$ is a state-feedback gain matrix for which the matrix $\bar{A} = A + BK$ is Schur (all eigenvalues are in the interior of the unit disk in the complex plane), $X(k)$ is the current spacecraft state, $X_n(k+\delta) \in \mathcal{N}$ is a time-varying reference along the NMT, and $\delta \in \mathbb{Z}$ is a shift which gives the controller set-point at the first time-instant the controller is switched to the specified NMT as the target.

## 2.4 Closed Loop Dynamics

Combining (1) and (9), the closed loop dynamics are given by

$$X(k+1) = \bar{A}X(k) - BKX_n(k+\delta), \tag{10}$$

where $\bar{A} = A + BK$. Defining the state error as $e(k, \delta) = X(k) - X_n(k+\delta)$, the error dynamics are given by

$$e(k+1, \delta) = \bar{A}e(k, \delta). \tag{11}$$

In the subsequent developments, the notation for the state error is simplified by omitting $\delta$, i.e., $e(k) = e(k, \delta)$, and $e(k+1) = e(k+1, \delta)$.

## 2.5 Constraints

Two constraints are considered. Firstly, the thrust is limited as

$$\|u(k)\|_\infty - u_{max} \leq 0, \tag{12}$$

where $\|\cdot\|_\infty$ denotes the infinity-norm. This constraint is equivalently stated as

$$\eta_i^T K(X(k) - X_n(k+\delta)) \leq u_{max}, \ i = 1, 2, ..., 6, \tag{13}$$

where $\eta_i$ are the vectors corresponding to the vertices of the unit infinity-norm hypercube, and $u_{max}$ is the norm bound specified by the mission designer.

Secondly, the satellite is required to stay out of one or more prescribed exclusion zones. These exclusion zones could, for example, represent the locations of other spacecraft or obstacles which must be avoided. The exclusion zones are modeled as ellipsoidal sets centered at specified points $s_i \in \mathbb{R}^3$. The $i$th exclusion zone is defined as

$$\mathcal{O}_i(s_i, S_i) = \{X \in \mathbb{R}^6 \mid (\Phi X - s_i)^T S_i (\Phi X - s_i) \leq 1\}, \tag{14}$$

where $S_i = S_i^T > 0$ is a shape matrix based on characteristics of the obstacle, including any uncertainty in its position, and the matrix $\Phi = [I_{3\times 3} \ 0_{3\times 3}]$ isolates the position from the state vector. The constraints on the spacecraft's position based on the $l$ exclusion zones are given by $X(k) \notin \mathcal{O}_i(s_i, S_i), i = 1, 2, ..., l$ which is equivalent to the inequality constraints,

$$1 - (\Phi X(k) - s_i)^T S_i (\Phi X(k) - s_i) \leq 0, \ i = 1, 2, ..., l. \tag{15}$$

Note that the constraints are given as inequalities in (12) and (15) to facilitate their use in the simulation results given in Section 5.

# 3 Safe Positively Invariant Tubes for Closed NMTs

In this section, safe, positively-invariant tubes are defined for closed NMTs. In this context, "safe" (constraint admissible) implies that constraints are satisfied point-wise within the tube, and "positively-invariant" implies that if the spacecraft state is within the tube at a given time instant, and the spacecraft motion is governed by the closed-loop dynamics (10), then it will remain within the tube for all future time. In the following subsections, this tube is formed by generating ellipsoidal sets about each point along the NMT $\mathcal{N}$, and then adjusting the sizes of these sets such that the tube formed by their union is both safe and positively-invariant.

## 3.1 Safe Sets

An ellipsoidal set, centered at the point $X_n(k) \in \mathcal{N}$ with scale factor $\rho_k \geq 0$, is defined as

$$\mathcal{E}_{k,\mathcal{N}} = \{X \in \mathbb{R}^6 \mid (X - X_n(k))^T P(X - X_n(k)) \leq \rho_k\}, \tag{16}$$

where the shape matrix $P = P^T > 0$ is chosen to satisfy the discrete Lyapunov inequality,

$$(A + BK)^T P(A + BK) - P < 0, \tag{17}$$

and where $A$ and $B$ are the discrete-time state and input matrices defined in (3) and (4), respectively, and $K$ is the state-feedback gain matrix defined in (9). With the same motivation as in [22], where safe sets are formed around forced equilibrium points, the set $\mathcal{E}_{k,\mathcal{N}}^s$ defined by (16) with $\rho_k = \rho_k^s$, i.e.,

$$\mathcal{E}_{k,\mathcal{N}}^s = \{X \in \mathbb{R}^6 \mid (X - X_n(k))^T P(X - X_n(k)) \leq \rho_k^s\}, \tag{18}$$

is safe if the scale factor $\rho_k^s$ is set to the largest possible value such that both:

a) the control constraint (12) is satisfied point-wise within the set with $\delta = 0$, i.e.,

$$\|u(k)\|_\infty = \|K(X - X_n(k))\|_\infty - u_{max} \leq 0 \text{ for all } X \in \mathcal{E}_{k,\mathcal{N}}^s \text{ and,}$$

b) the exclusion zone constraints (15) are satisfied point-wise within the set, i.e.,

$$1 - (\Phi X - s_i)^T S_i (\Phi X - s_i) \leq 0, \ i = 1, 2, ..., l \text{ for all } X \in \mathcal{E}_{k,\mathcal{N}}^s.$$

The scale factor $\rho_k^s$ is determined by first calculating the maximum possible scale factor for which item (a) holds, denoted by $\rho_u$, and for which item (b) separately holds, denoted by $\rho_{Oi,k}$. Then, $\rho_k^s$ is selected to be

$$\rho_k^s = \min\{\rho_u, \ \rho_{Oi,k}, \ i = 1, 2, ..., l\}. \tag{19}$$

### 3.1.1 Maximum Scale Factor Considering the Control Constraint

The control limit on the scale factor $\rho_k^s$, denoted $\rho_u$, is found by solving, for $i = 1, 2, ..., 6$, the following convex optimization problem:

$$\begin{aligned} \underset{X}{\text{maximize}} & \quad \eta_i^T K(X - X_n(k)) \\ \text{subject to} & \quad \tfrac{1}{2}(X - X_n(k))^T P(X - X_n(k)) \leq \alpha. \end{aligned} \tag{20}$$

If a value for $\alpha$ is found such that the solutions $X_i^*$ of (20) satisfy $\max_i\{\eta_i^T(K(X_i^* - X_n(k)))\} = u_{max}$, then $\rho_u = 2\alpha$.

The solution to (20) is obtained following the method developed in [22]. The matrix $P$ is diagonalized as $P = V^T \Lambda V$ where $V$ is orthogonal and $\Lambda$ is a diagonal matrix with eigenvalues of $P$ on the diagonal. Next, by defining parameters $\zeta_i$ and $h_i$ as

$$X - X_n(k) = V^T \Lambda^{-\frac{1}{2}} \zeta_i, \tag{21}$$

and

$$h_i^T = \eta_i^T K V^T \Lambda^{-\frac{1}{2}}, \tag{22}$$

the optimization problem (20) is re-formulated as

$$\begin{aligned} \underset{\zeta_i}{\text{maximize}} & \quad h_i^T \zeta_i \\ \text{subject to} & \quad \tfrac{1}{2} \zeta_i^T \zeta_i \leq \alpha, \end{aligned} \tag{23}$$

which has the solution of

$$\zeta_i = \frac{h_i}{\|h_i\|_2} \sqrt{2\alpha}. \tag{24}$$

Therefore, the control limit on the scale factor $\rho_k^s$ is given by

$$\rho_u = \min_i \frac{u_{max}^2}{\|h_i\|_2^2}. \qquad (25)$$

While the scale factor $\rho_u$ corresponds to the largest ellipsoidal set for which the control constraint is satisfied point-wise within the set, the control constraint will also be satisfied point-wise in any set with $\rho_k^s \leq \rho_u$.

**Remark 2:** Note that the scale factor $\rho_u$ is independent of the point $X_n(k)$. Hence, $\rho_u$ is constant for all $X_n(k)$ and must only be calculated once. This differs from the case considered in [22] where the controller set points were forced equilibria and the value of $\rho_u$ depended on the chosen forced equilibria.

### 3.1.2 Maximum Scale Factor Considering Exclusion Zone Constraints

The maximum size scale factor, $\rho_{Oi,k}$, considering the $i$th exclusion zone constraint is determined as the solution to a convex optimization problem in which the minimum sized ellipsoid, centered at $X_n(k)$, is sought which shares a common point with the exclusion zone, $\mathcal{O}_i(s_i, S_i)$. This is accomplished by solving

$$\begin{array}{rl} \underset{\rho_{Oi,k},\, X}{\text{minimize}} & \rho_{Oi,k} \\ \text{subject to} & (X - X_n(k))^T P (X - X_n(k)) \leq \rho_{Oi,k}, \\ & (\Phi X - s_i)^T S_i (\Phi X - s_i) \leq 1. \end{array} \qquad (26)$$

The solution to (26) is obtained via Karush-Kuhn-Tucker (KKT) conditions, [40], following the method used in [22].

Note that while the scale factor $\rho_{Oi,k}$ corresponds to the largest ellipsoidal set for which the $i$th exclusion zone constraint is satisfied point-wise within the set, the $i$th exclusion zone constraint will also be satisfied point-wise in any set with $\rho_k^s \leq \rho_{Oi,k}$.

**Remark 3:** If the point $X_n(k)$ lies within a keep-out zone, i.e., $X_n(k) \in \mathcal{O}_i(S_i, s_i)$, then no safe set may be formed. In this case $\rho_{Oi,k}$ is set to 0.

## 3.2 Safe Tubes

A safe tube centered on the NMT $\mathcal{N}$ is defined by

$$\mathcal{T}_\mathcal{N}^s = \bigcup_{k \in [0, k_{max}]} \mathcal{E}_{k,\mathcal{N}}^s. \qquad (27)$$

This tube is safe in the sense that for all $X(k) \in \mathcal{T}_\mathcal{N}^s$, the exclusion zone constraints (15) are satisfied and there exists $\delta \in \mathbb{Z}$ such that the control constraint (12) is satisfied.

Figure 2 shows three orthographic views of the projection of the 6-dimensional tube $\mathcal{T}_\mathcal{N}^s$ onto the position space for an example closed NMT. The tube $\mathcal{T}_\mathcal{N}^s$ was formed considering a control constraint of the form (12) and a single exclusion zone centered at the origin. In Figure 2, different colors correspond to different ellipsoids $\mathcal{E}_{k,\mathcal{N}}^s \subset \mathcal{T}_\mathcal{N}^s$.

If the spacecraft initial state $X(0) \in \mathcal{T}_\mathcal{N}^s$, then, with a suitable choice of $\delta$, constraints are guaranteed to be satisfied at that instant. However, there is no *a priori* guarantee that constraints will be satisfied for $k > 0$. The next subsection develops methods to construct a new tube, $\mathcal{T}_\mathcal{N}$, such that $\mathcal{T}_\mathcal{N}$ is both safe and positively-invariant, guaranteeing constraints will be satisfied for all future time.

## 3.3 Safe, Positively-Invariant Tubes

A safe, positively invariant tube $\mathcal{T}_\mathcal{N}$ is developed by generating a new set of scale factors $\rho_k$ from $\rho_k^s$ such that the property of positive-invariance holds and $\rho_k \leq \rho_k^s$ for all $k \in [0, k_{max}]$, hence $\mathcal{T}_\mathcal{N} \subset \mathcal{T}_\mathcal{N}^s$ and the safety of the tube is maintained.

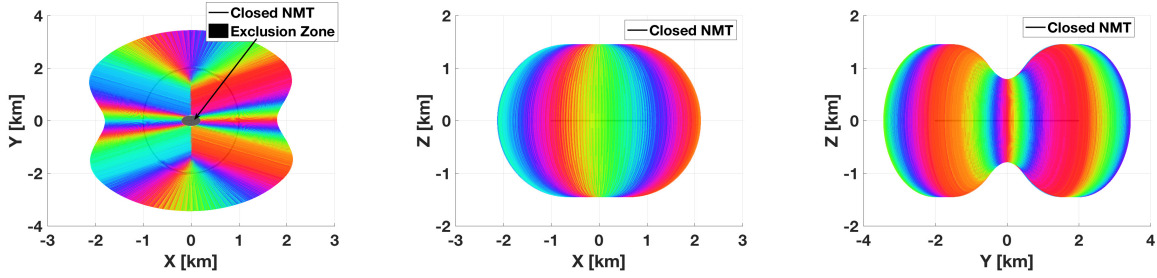

Fig. 2 Visualization of the safe tube $\mathcal{T}_\mathcal{N}^s$ projected onto $\mathbb{R}^3$

To guarantee all constraints are satisfied, it must hold that the spacecraft state vector is always within the ellipsoidal set corresponding to the current controller set point, i.e., $X \in \mathcal{E}_{k+\delta,\mathcal{N}}$. Therefore, the following definition for positive invariance is used:

**Definition 1:** Given a NMT $\mathcal{N}$, a tube

$$\mathcal{T}_\mathcal{N} = \bigcup_{k \in [0, k_{max}]} \mathcal{E}_{k,\mathcal{N}}, \tag{28}$$

is positively-invariant with respect to the closed loop dynamics given by (10) if there exists a $\delta \in \mathbb{Z}$ such that

$$X(k_1) \in \mathcal{E}_{k_1+\delta,\mathcal{N}} \implies X(k_2) \in \mathcal{E}_{k_2+\delta,\mathcal{N}} \ \forall \ k_2 \geq k_1, \ k_1, k_2 \in \mathbb{Z}_{\geq 0}. \tag{29}$$

**Remark 4:** Note that the definition of positive invariance in (29) implies $\mathcal{T}_\mathcal{N}$, as a set, is positively invariant with the appropriate selection of the control law (9) and $\delta$.

Two theorems are now presented which give conditions on the values of $\rho_k$ that result in a positively-invariant tube. These theorems can then be used to generate values for the scale factors $\rho_k$ from the safe scale factors $\rho_k^s$. The condition in Theorem 1 is a sufficient condition, and is conservative. The condition in Theorem 2 is both necessary and sufficient for positive invariance. However, applying Theorem 2 to determine $\rho_k$ requires slightly more computation time compared to Theorem 1.

### 3.3.1 Conditions for Positive Invariance

Two assumptions are needed for both Theorems 1 and 2. These assumptions were previously introduced in Section 2 and are re-stated here for clarity:

A1) The closed loop dynamics and error dynamics are given by (10) and (11), respectively.

A2) The ellipsoidal set shape matrix $P = P^T > 0$ is chosen such that $\bar{A}^T P \bar{A} - P = -Q$, $Q = Q^T > 0$.

The condition for positive invariance in Theorem 1 is developed by leveraging Assumption (A2) which ensures that the state error $e(k)$ always decays to successively smaller ellipsoids as time progresses, i.e., $e(k)^T P e(k) - e(k+1)^T P e(k+1) > 0$ for all $k \in \mathbb{Z}_{\geq 0}$.

**Theorem 1:** Suppose Assumptions (A1) and (A2) hold. Then the tube $\mathcal{T}_\mathcal{N}$ is positively-invariant if

$$\rho_{k_1} \leq \rho_{k_2} \text{ whenever } k_1 \leq k_2. \tag{30}$$

**Proof:** Without loss of generality, let $k_2 = k_1 + 1$ and let $\delta = 0$. Assume $X(k_1) \in \mathcal{E}_{k_1,\mathcal{N}}$, therefore $e(k_1)^T P e(k_1) \leq \rho_{k_1}$. By Assumptions (A1) and (A2),

$$e(k_2)^T P e(k_2) - e(k_1)^T P e(k_1) = e(k_1)^T (\bar{A}^T P \bar{A} - P) e(k_1) = -e(k_1)^T Q e(k_1). \tag{31}$$

Since $Q = Q^T > 0$,

$$e(k_2)^T P e(k_2) - e(k_1)^T P e(k_1) \leq -\lambda_{min}(Q) \|e(k_1)\|_2^2, \tag{32}$$

where $\lambda_{min}(Q) \in \mathbb{R}_{>0}$ is the minimum eigenvalue of the matrix Q. Therefore, $e(k_2)^T P e(k_2) < e(k_1)^T P e(k_1)$. If $\rho_{k1} \leq \rho_{k2}$, it is guaranteed that $e(k_2)^T P e(k_2) < \rho_{k2}$. Therefore, $X(k_2) \in \mathcal{E}_{k_2, \mathcal{N}}$. □

Theorem 1 may be applied to generate scale factors $\rho_k$ from $\rho_k^s$ such that the resulting tube is both safe and positively-invariant. It is shown later in Subsection 3.3.2 that the calculations required are minimal; however, this tube may be much smaller than the initial safe tube, limiting its utility in trajectory planning. To address this limitation, Theorem 2 is developed which gives a necessary and sufficient condition for positive invariance. Hence, a tube generated using Theorem 2 is as large as possible given an ellipsoidal shape matrix $P$.

The idea behind the condition for positive invariance given in Theorem 2 is to determine the smallest possible amount by which the error ellipsoid $e^T P e$ will shrink over one discrete-time step. Then, by ensuring that the ellipsoidal sets, with size defined by scale factors $\rho_k$, along the NMT do not shrink by more than this amount, the resulting tube is guaranteed to be positively-invariant. This idea is first applied in Lemma 1 below. The proof of Lemma 1 relies on the following proposition.

**Proposition 1:** Let $D \subset \mathbb{R}^n$ be a compact, convex set with a non-empty interior, and let $f(X) : D \to \mathbb{R}$ be convex. If $f(X) \leq a$ when $X \in \partial D$, i.e., $X$ is on the boundary of $D$, then $f(X) \leq a$ for all $X \in D$.

**Proof:** Follows from the standard fact that the maximum of a convex function over a compact, convex set occurs on the boundary of the set, see Corollary 32.2.1 in [41].

**Lemma 1:** Suppose Assumptions (A1) and (A2) hold. The tube $\mathcal{T}_\mathcal{N}$ is positively-invariant if and only if

$$\rho_{k+1} \geq \rho_k - \bar{d}(\rho_k) \ \forall \ k \in \mathbb{Z}_{\geq 0}, \tag{33}$$

where

$$\bar{d}(\rho_k) = \min_{e(k)} \quad e(k)^T P e(k) - e(k+1)^T P e(k+1) \tag{34}$$
$$\text{subject to} \quad e(k)^T P e(k) = \rho_k.$$

**Proof:** Without loss of generality, let $\delta = 0$. To prove sufficiency, suppose (33), (34) hold. Assume $X(k) \in \mathcal{E}_{k,\mathcal{N}}$, therefore $e(k)^T P e(k) \leq \rho_k$. Define $D = \{e \mid e^T P e \leq \rho_k\}$. The set $D$ is a closed and bounded subset of $\mathbb{R}^6$, i.e., compact. It is also convex since it is a sub-level set of the strictly convex function $e^T P e$, $P > 0$. Define $f(e(k)) = e(k)^T \bar{A}^T P \bar{A} e(k) = e(k+1)^T P e(k+1)$. The matrix $\bar{A}$ is invertible, hence $\bar{A}^T P \bar{A} > 0$, and $f(e(k))$ is convex. From (34), $f(e(k)) \leq \rho_k - \bar{d}(\rho_k)$ for all $e(k) \in \partial D$. Therefore, by Proposition 1,

$$f(e(k)) = e(k+1)^T P e(k+1) \leq \rho_k - \bar{d}(\rho_k), \tag{35}$$

for all $e(k) \in D$. Combining (33) and (35) yields

$$e(k+1)^T P e(k+1) \leq \rho_k - \bar{d}(\rho_k) \leq \rho_{k+1}. \tag{36}$$

Therefore, $X(k+1) \in \mathcal{E}_{k+1,\mathcal{N}}$. The same arguments may be repeated to show $X(k+2) \in \mathcal{E}_{k+2,\mathcal{N}}$, ... , $X(k+n) \in \mathcal{E}_{k+n,\mathcal{N}}$, therefore $\mathcal{T}_\mathcal{N}$ is positively-invariant.

The proof for necessity is by contradiction. Suppose $\mathcal{T}_\mathcal{N}$, defined by (28), (16), is positively-invariant, but (33) does not hold, i.e.,

$$\rho_{k+1} < \rho_k - \bar{d}(\rho_k). \tag{37}$$

Consider $X(k)$ such that $e(k)^T P e(k) = \rho_k$. Then, (34) gives $e(k+1)^T P e(k+1) \leq \rho_k - \bar{d}(\rho_k)$. Hence, there exists a $X(k)$ such that $X(k) \in \mathcal{E}_{k,\mathcal{N}}$ and $e(k+1)^T P e(k+1) = \rho_k - \bar{d}(\rho_k) > \rho_{k+1}$, hence $X(k+1) \notin \mathcal{E}_{k+1,\mathcal{N}}$, contradicting the assumption that $\mathcal{T}_\mathcal{N}$ is positively invariant. □

Note that while Lemma 1 provides a necessary and sufficient condition for positive invariance, using the condition given in (33) to generate $\rho_k$ from $\rho_k^s$ may result in a $\rho_k > \rho_k^s$. This is because to apply condition (33), the parameter $\bar{d}(\rho_k)$ must be calculated from a given value of $\rho_k$, and then used to potentially adjust $\rho_k^s$ upward. Therefore, while the tube $\mathcal{T}_\mathcal{N}$ is guaranteed to be positively-invariant, it may not be safe. To generate a positively-invariant tube which is also guaranteed to be safe, an equivalent condition to (33) is derived by calculating the minimum change in $\rho_k$ over one discrete-time step by looking backward in time. Before stating Theorem 2, the following Lemma is presented to establish the equivalence of the conditions given in Lemma 1 and Theorem 2.

**Lemma 2:** Suppose Assumptions (A1) and (A2) hold. Then,

$$\rho_k \leq \rho_{k+1} + d(\rho_{k+1}) \ \forall \ k \in \mathbb{Z}_{\geq 0}, \tag{38}$$

where

$$d(\rho_{k+1}) = \min_{e(k+1)} e(k)^T P e(k) - e(k+1)^T P e(k+1) \\ \text{subject to } e(k+1)^T P e(k+1) = \rho_{k+1}, \tag{39}$$

if and only if

$$\rho_{k+1} \geq \rho_k - \bar{d}(\rho_k) \ \forall \ k \in \mathbb{Z}_{\geq 0}, \tag{40}$$

where

$$\bar{d}(\rho_k) = \min_{e(k)} e(k)^T P e(k) - e(k+1)^T P e(k+1) \\ \text{subject to } e(k)^T P e(k) = \rho_k. \tag{41}$$

**Proof:** See Appendix B.1.

**Theorem 2:** Suppose Assumptions (A1) and (A2) hold. The tube $\mathcal{T}_\mathcal{N}$ is positively-invariant if and only if

$$\rho_k \leq \rho_{k+1} + d(\rho_{k+1}) \ \forall \ k \in \mathbb{Z}_{\geq 0}, \tag{42}$$

where

$$d(\rho_{k+1}) = \min_{e(k+1)} e(k)^T P e(k) - e(k+1)^T P e(k+1) \\ \text{subject to } e(k+1)^T P e(k+1) = \rho_{k+1}. \tag{43}$$

**Proof:** The proof follows directly from Lemmas 1 and 2. □

### 3.3.2 Application of Theorems 1 and 2

Both Theorems 1 and 2 can be used to generate safe, positively-invariant tubes by determining values for scale factors $\rho_k$ from the safe scale factors $\rho_k^s$. Below, procedures are given to accomplish this.

**Application of Theorem 1**

Due to the periodicity of the NMTs, in order to satisfy the conditions of Theorem 1, a safe, positively-invariant tube $\mathcal{T}_\mathcal{N}$ may be formed using the following procedure.
**Procedure 1:**

1. Set $\rho_k = \min_{k \in [0, k_{max}]} \rho_k^s$ for $k \in [0, k_{max}]$.

It is clear that defining $\mathcal{T}_\mathcal{N}$ using this method requires minimal calculations, however, doing so may result in a tube much smaller than the safe tube $\mathcal{T}_\mathcal{N}^s$. This will occur, for example, for an NMT which passes very near an exclusion zone.

**Application of Theorem 2**

By applying Theorem 2, a larger tube may be formed. In order to do so, the Quadratically Constrained Quadratic Program (QCQP) given in (43) must be solved. This is done efficiently by first converting the QCQP to a Linear Program (LP), see Appendix C. Then, the scale factors $\rho_k$ are obtained from $\rho_k^s$ using the following procedure.
**Procedure 2:**

1. Start at $k$ such that

$$k + 1 = \underset{k \in [0, k_{max}]}{\operatorname{argmin}} (\rho_k^s), \tag{44}$$

and set $\rho_{k+1} = \rho_{k+1}^s$. Note that if multiple $k$ satisfy (44), any such $k$ can be chosen as the starting location.

2. Determine $d(\rho_{k+1})$ and set

$$\rho_k = \begin{cases} \rho_k = \rho_k^s & \text{if } \rho_k^s \leq \rho_{k+1} + d(\rho_{k+1}), \\ \rho_k = \rho_{k+1} + d(\rho_{k+1}) & \text{otherwise.} \end{cases} \quad (45)$$

3. Increment $k = k - 1$ and repeat step 2. When $k = -1$, set $k = k_{max}$ and continue until returning to the starting index.

**Remark 5:** Procedure 2 yields the largest safe, positively-invariant tube comprised of ellipsoidal sets with shape matrix $P$. Specifically, setting $\rho_k$ using Procedure 2 results in the largest possible $\rho_k$ such that both $\rho_k \leq \rho_k^s$ and condition (42) from Theorem 2 holds.

**Remark 6:** Note that if any $\rho_k^s = 0$, then the NMT passes through an exclusion zone and motion along the NMT is not safe. In this case, using either Procedure 1 or 2 to select $\rho_k$ will result in $\rho_k = 0$ for all $k \in [0, k_{max}]$. This is desirable as it ensures that the NMT will not be included in any trajectories planned using the methods described in Section 4.

Using Procedures 1 or 2, a safe, positively-invariant tube with $\rho_k > 0$ for all $k \in [0, k_{max}]$ can be formed about any trajectory which does not enter any exclusion zones. This statement is presented formally in Theorems 3 and 4. In these theorems, the notation $\emptyset$ is used to denote the empty set.

**Theorem 3:** Suppose Assumptions (A1) and (A2) hold, $u_{max} > 0$, and $\mathcal{N} \cap \mathcal{O}_i(s_i, S_i) = \emptyset$, for $i = 1, 2, ..., l$. Then, there exist $\rho_k > 0$ for all $k \in [0, k_{max}]$ obtained using Procedure 1 such that $\mathcal{T}_\mathcal{N}$ is safe and positively-invariant.
**Proof:** See Appendix B.2

**Theorem 4:** Suppose Assumptions (A1) and (A2) hold, $u_{max} > 0$, and $\mathcal{N} \cap \mathcal{O}_i(s_i, S_i) = \emptyset$ for $i = 1, 2, ..., l$. Then, there exist $\rho_k > 0$ for all $k \in [0, k_{max}]$ obtained using Procedure 2 such that $\mathcal{T}_\mathcal{N}$ is safe and positively-invariant.
**Proof:** See Appendix B.3

### 3.3.3 Example showing implementation of Procedures 1 and 2

The safe, positively-invariant tubes for an example NMT generated using Procedures 1 and 2 are illustrated in Figure 3. Figure 3a contains plots of both $\rho_k^s$ and $\rho_k$ showing how the values of $\rho_k^s$ are adjusted using Procedures 1 and 2. Note that the values for $\rho_k$ generated using Procedure 1 are all constant and equal to the minimum value of $\rho_k^s$. It is clear that using Procedure 1, the values of $\rho_k$ are limited by the minimum value of $\rho_k^s$.

The values of $\rho_k$ generated using Procedure 2 are also less than or equal to the corresponding $\rho_k^s$ value, however they increase and decrease along the trajectory such that the less conservative condition of Theorem 2 is satisfied. The starting location for Procedure 2 is denoted in Figure 3a by a $\square$. Because Procedure 2 starts at the minimum value of $\rho_k^s$ and proceeds backward in terms of the discrete-time instances, and because the maximum value for $\rho_k$ is limited by the value of $\rho_{k+1}$, scale factors $\rho_k$ determined by Procedure 2 are also limited by the minimum value of $\rho_k^s$.

Figures 3b, 3c and 3d show projections of the safe tube $\mathcal{T}_\mathcal{N}^s$, and safe, positively-invariant tubes $\mathcal{T}_\mathcal{N}^1$ and $\mathcal{T}_\mathcal{N}^2$ generated using Procedures 1 and 2, respectively. Note that $\mathcal{T}_\mathcal{N}^1 \subset \mathcal{T}_\mathcal{N}^2 \subset \mathcal{T}_\mathcal{N}^s$. Additionally, the tube $\mathcal{T}_\mathcal{N}^2$ shown in Figure 3b is very similar to the safe tube $\mathcal{T}_\mathcal{N}^s$ shown in Figure 3d, illustrating that Procedure 2 makes relatively small adjustments to the scale factors $\rho_k^s$.

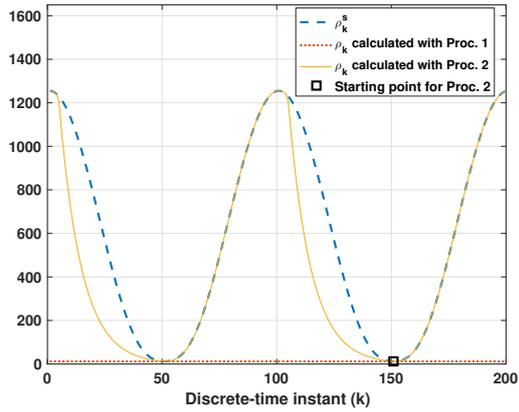

(a) Plot of $\rho_k^s$, and $\rho_k$ obtained by adjustments to $\rho_k^s$ using Procedures 1 and 2

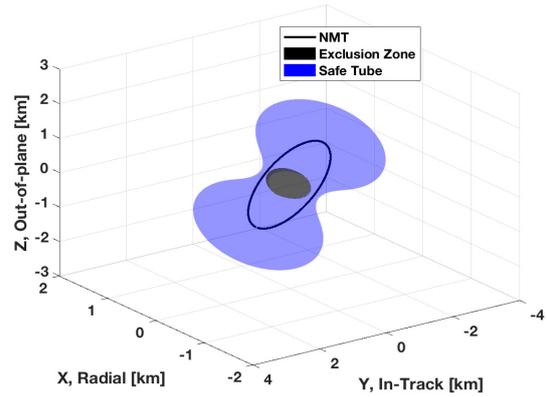

(b) Projection of safe tube corresponding to $\rho_k^s$

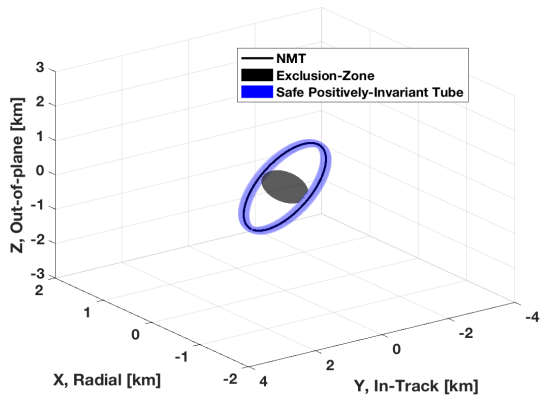

(c) Projection of safe, positively-invariant tube corresponding to $\rho_k$ calculated with Procedure 1

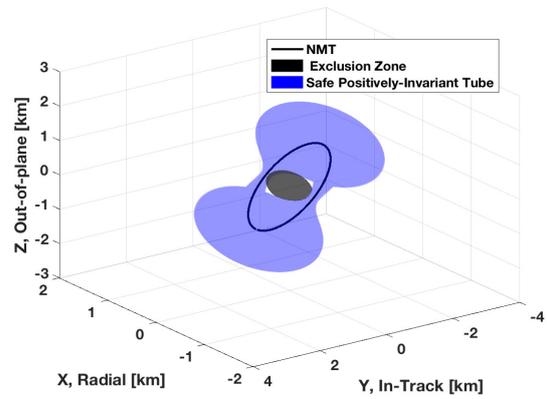

(d) Projection of safe, positively-invariant tube corresponding to $\rho_k$ calculated with Procedure 2

Fig. 3 Example illustrating the application of Procedures 1 and 2

# 4 Virtual Net for Safe Trajectory Planning

Using the safe, positively-invariant tubes defined above, safe trajectories are planned between desired closed NMTs. Given a set of $m$ closed NMTs, $\mathcal{M} = \{\mathcal{N}_1, \mathcal{N}_2, ..., \mathcal{N}_m\}$, a virtual net is formed. This virtual net consists of a directed graph with one node corresponding to each closed NMT, $\mathcal{N}_i \in \mathcal{M}$, and is represented by an adjacency matrix and connection array, defined below. A single node is sufficient for each NMT due to the periodicity, i.e., if the spacecraft is able to reach a single point along a closed NMT, it is able to reach all points along the closed NMT given sufficient time. Using the safe, positively-invariant tube for each node, the adjacency of nodes in the virtual net is determined, along with connection information which consists of both the starting point and the initial controller reference point for a transfer between NMTs. Note that the adjacency information only determines if a safe transfer is possible from one node to another, while the connection information provides the starting location and controller set point to actually execute a transfer. After the adjacency and connection information has been determined, safe trajectories are generated using efficient graph search algorithms such as Dijkstra's algorithm [17].

## 4.1 Virtual Net Adjacency

The definition of adjacency of two nodes $\mathcal{N}_1$ and $\mathcal{N}_2$ is given as:
**Definition 2:** $\mathcal{N}_1$ is adjacent to $\mathcal{N}_2$ if there exists $k_1 \leq k_{max}$ and $\hat{k}_2 \leq k_{max}$ such that

$$\mathcal{B}(X_{n1}(k_1), \gamma_1) \subset \mathcal{E}_{\hat{k}_2, \mathcal{N}_2}, \tag{46}$$

where $\mathcal{B}(Z, \gamma) = \{X \mid \|X - Z\|_2 \leq \gamma\}$, and $\gamma_1$ is a small parameter chosen by the mission designer. Note that choosing a larger value for $\gamma_1$ may result in fewer pairs of adjacent nodes. See Remarks 7 and 8 for additional discussion regarding the choice of $\gamma_1$. Figure 4 shows a sketch illustrating the parameters used in the adjacency definition (46). The requirement $\mathcal{B}(X_{n1}(k_1), \gamma_1) \subset \mathcal{E}_{\hat{k}_2, \mathcal{N}_2}$ in (46) along with the positive-invariance of $\mathcal{T}_{\mathcal{N}_2}$ ensures that a transfer from $\mathcal{N}_1$ to $\mathcal{N}_2$ may be executed without violating constraints by setting the controller reference point to $X_{n2}(\hat{k}_2) \in \mathcal{N}_2$ when the spacecraft is near $X_{n1}(k_1) \in \mathcal{N}_1$.

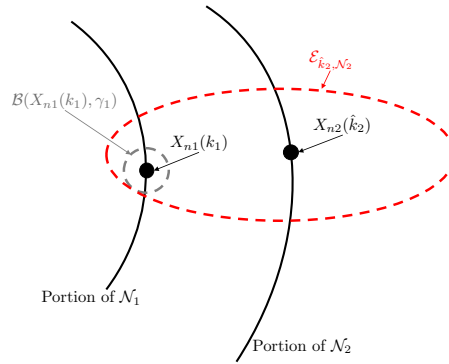

**Fig. 4** Illustration of the parameters used in adjacency definition (46)

**Remark 7:** It is also possible to define adjacency between nodes by replacing the requirement $\mathcal{B}(X_{n1}(k_1), \gamma_1) \subset \mathcal{E}_{\hat{k}_2, \mathcal{N}_2}$ with

$$X_{n1}(k_1) \in \text{Int}(\mathcal{E}_{\hat{k}_2, \mathcal{N}_2}), \tag{47}$$

which is consistent with the definition of adjacency between forced equilibrium points in [22]. If a pair of nodes are adjacent by (47), then there exists a $\gamma_1 > 0$ such that the pair is also adjacent by definition (46). The adjacency definition (46) is used here because this allows for more control in defining the switching behavior, i.e., the criteria used to determine when the controller reference point is switched to the next node, and ensures better match of predicted fuel consumption in the graph search optimization and actual fuel consumption. Note also that in (46), $\mathcal{B}(X_{n1}(k_1), \gamma_1)$ may be replaced with any bounded set containing $X_{n1}(k_1)$ in its interior. Finally, in (46), the parameter $\gamma_1$ may be chosen to be 0. This choice may be made to simplify adjacency calculations at the expense of possible constraint violation since the spacecraft only asymptotically approaches the NMT under the control law (9). Simulations show that this simplified implementation rarely leads to constraint violation.

### 4.1.1 Adjacency Matrices

An unweighted adjacency matrix, $\Pi$, is generated using the adjacency definition (46) as follows. To determine if $\mathcal{N}_i$ is adjacent to $\mathcal{N}_j$ for $i, j \in [1, 2, ..., m]$, a grid search over all $X_{ni}(k_i) \in \mathcal{N}_i$ and $X_{nj}(\hat{k}_j) \in \mathcal{N}_j$ is performed until the first pair satisfying (46) is found. If $\mathcal{N}_i$ is adjacent to $\mathcal{N}_j$ for $i, j \in [1, 2, ..., m]$, the corresponding matrix element is set to 1, i.e., $\Pi(i, j) = 1$. If $\mathcal{N}_i$ is not adjacent to $\mathcal{N}_j$ for $i, j \in [1, 2, ..., m]$, the corresponding matrix element is set to $+\infty$. Note that if $\mathcal{N}_i$ is adjacent to $\mathcal{N}_j$, it does not imply, in turn, that $\mathcal{N}_j$ is adjacent to $\mathcal{N}_i$.

The unweighted adjacency matrix can be used to generate safe trajectories between NMTs, however, these trajectories may not be fuel efficient. To generate trajectories with decreased fuel consumption, i.e., less control usage, a weighted adjacency matrix, $\Pi_w$, is generated by determining the most control-efficient transfer between connected NMTs as follows. A grid search over all $X_{ni}(k_i) \in \mathcal{N}_i$ and $X_{nj}(\hat{k}_j) \in \mathcal{N}_j$ is performed. For each pair of state vectors $X_{ni}(k_i)$, $X_{nj}(\hat{k}_j)$ satisfying the adjacency criteria given in (46), a transfer trajectory is calculated using (1), (9), starting at initial point $X(0) = X_{ni}(k_i)$ and initial controller reference $X_n(0) = X_{nj}(\hat{k}_j)$. The trajectory is propagated until the current state is within a small neighborhood of the controller reference, i.e., until $k = \bar{k}$ such that $X(\bar{k}) \in \mathcal{B}(X_{nj}(\bar{k} + \hat{k}_j), \gamma_2)$, where $\gamma_2$ is a small positive value chosen by the mission designer (see Remark 8). The cost of transition is calculated as the total control used over this trajectory, scaled by $\frac{1}{\Delta T}$,

$$u_{tot} = \sum_{k=0}^{\bar{k}} \|u(k)\|_1, \qquad (48)$$

where $\|\cdot\|_1$ denotes the 1-norm (note that the 1-norm is used here to represent fuel usage on a spacecraft with three thrusters (or thruster pairs) mounted orthogonal to each other). The weighted adjacency matrix $\Pi_w$ is formed by storing the lowest cost of transition between adjacent NMTs in the corresponding matrix element and, for NMTs which are not adjacent, the corresponding matrix element is set to $+\infty$.

## 4.2 Virtual Net Connection Information

To aid in trajectory planning, in addition to the weighted and unweighted adjacency matrices, weighted and unweighted connection arrays, $\Xi$ and $\Xi_w$, are formed. These connection arrays store the transfer starting location and the initial controller reference point used to execute the transfer between each pair of adjacent NMTs. Specifically, for each pair of adjacent NMTs, the corresponding element of the connection array consists of a vector containing the indices of the initial transfer location and initial controller reference point to be used to execute the transfer, i.e., $\Xi(i, j) = [k_i \ \hat{k}_j]$ where $X_{ni}(k_i)$ and $X_{nj}(\hat{k}_j)$ satisfy the adjacency criteria given in (46).

The parameters $k_i$ and $\hat{k}_j$ are determined for the unweighted and weighted connection arrays as follows. For the unweighted connection array, $\Xi$, $k_i$ and $\hat{k}_j$ can be any indices satisfying the adjacency criteria of (46), and are selected by performing a grid search over all $X_{ni}(k_i) \in \mathcal{N}_i$ and $X_{nj}(\hat{k}_j) \in \mathcal{N}_j$ until the first pair satisfying (46) is found. Note that a transfer executed using the information in the unweighted connection array is not guaranteed to be fuel efficient. For the weighted connection array, $\Xi_w$, a grid search over all $X_{ni}(k_i) \in \mathcal{N}_i$ and $X_{nj}(\hat{k}_j) \in \mathcal{N}_j$ is performed and, for each pair $X_{ni}(k_i)$ and $X_{nj}(\hat{k}_j)$ satisfying (46), the control-cost (48) is calculated. The values for $k_i$ and $\hat{k}_j$ which correspond to the lowest control-cost are stored in the weighted connection array, $\Xi_w$. Hence, a transfer executed using the information in the weighted connection array is expected to be the lowest control-cost trajectory between the specified NMTs $\mathcal{N}_i$ and $\mathcal{N}_j$ when the possible starting locations and controller set-points are confined to the sets $\mathcal{N}_i$ and $\mathcal{N}_j$, respectively, defined in (7).

## 4.3 Trajectory Planning

After forming the virtual net, safe trajectories can be planned on-line using the initial conditions for each NMT in the virtual net, and the unweighted or weighted adjacency matrix and connection array. Dijkstra's algorithm is applied to generate a sequence of nodes (NMTs) which connects given starting and ending NMTs. Note that Dijkstra's algorithm checks the adjacency of the entire virtual net to generate the sequence of nodes, and that the algorithm is complete, i.e., if a solution exists, the algorithm will return the solution.

After a sequence of NMTs has been obtained, the appropriate connection array is used to generate a safe trajectory by switching the controller reference to the next NMT in the sequence once the spacecraft reaches a small neighborhood of each transfer location. Specifically, for a spacecraft travelling toward NMT $\mathcal{N}_i$ before transferring to NMT $\mathcal{N}_j$, the connection array element $\Xi(i,j) = [k_i \ \hat{k}_j]$. The controller reference is switched to $X_{nj}(\hat{k}_j)$ at the first time-instant $\bar{k}$ when the spacecraft state vector satisfies

$$X(\bar{k}) \in \mathcal{B}(X_{ni}(k_i), \gamma_3), \tag{49}$$

where $\gamma_3$ is a small parameter chosen by the mission designer.

**Remark 8:** Note that in (49), $\mathcal{B}(X_{ni}(k_i), \gamma_3)$ may be replaced with any convex set containing $X_{ni}(k_i)$ in its interior. Choosing the sets used to define adjacency in (46) and switching in (49) to be the same, i.e, choosing $\gamma_1 = \gamma_3$, ensures all transfer trajectories between NMTs will satisfy constraints, and therefore that any trajectory consisting of multiple transfers between successive pairs of NMTs will be safe. Additionally, choosing $\gamma_1 = \gamma_2 = \gamma_3$ yields the best match between predicted fuel usage and actual fuel usage. However, as noted earlier, $\gamma_1$ may be chosen to be 0 to simplify adjacency calculations. Simulations show that choosing $\gamma_1 = 0$ and $\gamma_2 = \gamma_3$ to be small, but non-zero, rarely leads to constraint violations and provides a good estimate of fuel usage.

## 5 Simulations

Simulation case studies are now considered. Table 2 lists spacecraft parameters, nominal circular orbit parameters, constraints, and parameters used to determine adjacency, transfer costs and controller switching times. The state-feedback gain matrix $K$ for the controller (9) is an LQ gain matrix corresponding to the selection of state and control weighting matrices given by $Q_{LQ} = 100 \ \text{diag}(1, \ 1, \ 1, \ 1 \times 10^5, \ 1 \times 10^5, \ 1 \times 10^5)$ and $R_{LQ} = 2 \times 10^7 I_{3 \times 3}$. The shape matrix $P$ for the ellipsoidal set computations was chosen to be the solution to the discrete-time Riccati equation in the LQ problem [42]. The projection and visualization of ellipsoidal sets was accomplished using the Ellipsoidal Toolbox for Matlab® [43], and Dijkstra's algorithm was implemented using the MatlabBGL toolbox [44].

Table 2  Parameters used in simulations

| Parameter | Symbol | Value |
| --- | --- | --- |
| Spacecraft mass | $m$ | 140 kg |
| Nominal orbital radius for CW dynamics | $R_0$ | 7728 km |
| Mean motion | $\omega$ | 0.001027 rad/sec |
| Discrete-time update period | $\Delta T$ | 30.58 sec |
| Discrete-time index of final point (before repeating) on closed NMTs | $k_{max}$ | 199 |
| Maximum allowable control | $u_{max}$ | 0.005 kg·km/sec$^2$ (5 N) |
| Center of exclusion zone 1 | $s_1$ | $[0 \ 1 \ 0]^T$ km |
| Center of exclusion zone 2 | $s_2$ | $[0 \ -1 \ 0]^T$ km |
| Shape matrix for exclusion zones | $S_i, i = 1, 2$ | $\frac{1}{0.2^2} I_{3 \times 3}$ |
| Parameter used to determine adjacency | $\gamma_1$ | 0 |
| Parameter used to determine cost of transition between nodes | $\gamma_2$ | 0.0001 |
| Parameter used to determine controller switching times | $\gamma_3$ | 0.0001 |

## 5.1 Simulation Virtual Net

A set of 84 closed NMTs is used including 54 elliptical NMTs centered at the origin, 15 straight-line segment periodic NMTs and 15 stationary-point NMTs (in-track equilibria). These NMTs are chosen to be evenly spaced within a box of $3.5 \times 7 \times 10$ km in the $X$, $Y$, and $Z$ directions, respectively, centered at the origin. The 54 elliptical NMTs are chosen to have initial conditions corresponding to all combinations of parameters $b$, $\theta_1$ and $\theta_2$, defined in Appendix A, given by $b = \{0.5,\ 0.75,\ ...\ ,1.75\}$, $\theta_1 = \{45^o,\ 90^o,\ 135^o\}$, and $\theta_2 = \{-45^o,\ 0^o,\ 45^o\}$. The 15 straight line segment NMTs and 15 stationary point NMTs are chosen to be evenly spaced along the Y-axis with intersections at $y = \{-3.5,\ -3,\ ...\ ,3.5\}$. These 84 NMTs are shown in Figure 5. Two virtual nets are formed, corresponding to safe, positively-invariant tubes generated using

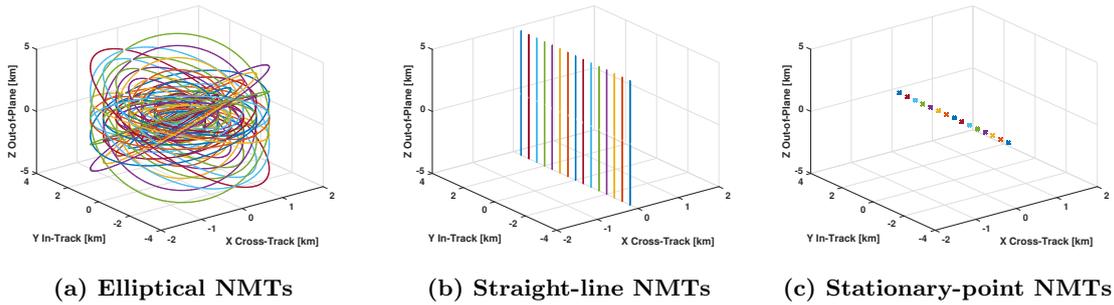

(a) Elliptical NMTs  (b) Straight-line NMTs  (c) Stationary-point NMTs

**Fig. 5** NMTs in the virtual net

Procedures 1 and 2, respectively.

**Remark 9:** There is much flexibility in choosing what types of NMTs, and how many NMTs, are included in the virtual net. Increasing the number of NMTs in the virtual net may increase the number of pairs of adjacent NMTs, increase the number of feasible trajectories between NMTs, and/or decrease the cost of trajectories connecting NMTs. These potential benefits come at the expense of additional computations required to generate the adjacency matrix and connection array. Development of methods to optimally chose NMTs for the virtual net will be considered in future work.

## 5.2 Description of Simulation Figures

In the results that follow, in addition to the trajectories and constraints, the parameter

$$w(k) = e(k)^T P e(k) - \rho_{k+\delta_i}^i \qquad (50)$$

is plotted, where $\rho_{k+\delta_i}^i$ is the ellipsoidal scale factor for the current controller set-point. This is done to demonstrate that the trajectories stay within the tubes $\mathcal{T}_{\mathcal{N}_i}$ at all times. Note that if $w(k) \leq 0$, then the current state $X(k)$ is within the safe-positively invariant tube corresponding to the current controller set-point, i.e., $X(k) \in \mathcal{T}_{\mathcal{N}_i}$.

For visual clarity, legends are not included in plots showing trajectories. In these plots, exclusions zones are shown by gray ellipsoids. The spacecraft trajectory is depicted as a solid pink line, the initial NMT is a dashed green line, the final NMT is a dashed red line and intermediate NMTs are depicted as dashed black lines. Transfer locations and initial controller reference points are depicted with * and □, respectively. The initial spacecraft position is shown as a green X, and the final position is shown with a red O.

## 5.3 Simulation Results

Figures 6a and 6b show a simulation where the trajectory is constructed using the virtual net calculated using Procedure 1. In this example, a safe trajectory is planned between two elliptical NMTs centered at the origin. The initial condition used to generate the initial and final NMTs are given by $\bar{X}_0^i = [0\ 1\ -1\ 0.0005\ 0\ -0.0007]^T$ and $\bar{X}_0^f = [0\ 3\ 0\ 0.0015\ 0\ 0.0015]^T$, respectively, with units of km for position and km/sec for velocity.

The trajectory calculated using the unweighted adjacency matrix and connection array, along with the corresponding constraints and $w(k)$, are shown in Figure 6a, while the same parameters calculated using the weighted adjacency matrix and connection array are shown in Figure 6b. Note that in each case, constraints are satisfied, i.e., all constraint values are $\leq 0$, and $w(k) \leq 0$ for the entire trajectory. Using the same initial and final NMTs, a simulation is run using a virtual net calculated using Procedure 2 and a weighted adjacency matrix, see Figure 6c. Note that this trajectory utilizes more intermediate NMTs then the trajectories calculated on the Procedure 1 virtual net and that the trajectory forms a spiral between successively larger NMTs for much of the trajectory. The total cost, $J$, of each trajectory, corresponding to the total control usage along that trajectory, i.e., $J = \Delta T \sum_0^{k_{final}} \|u(k)\|_1$, is shown in Table 3 with units of N· sec. As expected, the trajectories planned using the weighted adjacency matrices have lower total control costs than the trajectory planned using the unweighted matrices. The total cost for the trajectory planned using Procedure 2 and a weighted adjacency matrix is lower than either cost obtained using Procedure 1. This reduction in cost may be due to the increased adjacency of the virtual net compared to the virtual net calculated using Procedure 1. Specifically, in the virtual net formed using Procedure 1, there are 1501 pairs of adjacent nodes, while in the virtual net formed using Procedure 2, there are 2457 pairs of adjacent nodes.

Table 3  Cost comparison

| Virtual net | Adjacency and Connections | Shown in Figure | Cost |
|---|---|---|---|
| Theorem 1 | Unweighted | 6a | 1480 N·sec |
| Theorem 1 | Weighted | 6b | 951 N·sec |
| Theorem 2 | Weighted | 6c | 930 N·sec |

An additional advantage of the increased adjacency provided by the virtual net calculated using Procedure 2 is that safe trajectories may be planned between NMTs that are not possible using a virtual net calculated using Procedure 1. This is illustrated in Figures 7 and 8. Consider an NMT, denoted $\mathcal{N}_f$, that passes nearby an exclusion zone. Such a straight line segment NMT is plotted in Figure 7 as a solid red line. The corresponding safe, positively-invariant tube is plotted in purple. The tube for $\mathcal{N}_f$ generated using Procedure 1 is shown in Figure 7a. Note that the tube is small, and no nearby NMTs (in terms of distance, plotted as black dashed lines) pass through it. Hence, no other NMTs are adjacent to $\mathcal{N}_f$, and no trajectories ending on $\mathcal{N}_f$ are possible. Figure 7b shows that the tube for $\mathcal{N}_f$ generated using Procedure 2 is much larger, and connections to $\mathcal{N}_f$ are possible.

Figure 8 shows an example trajectory from a stationary point NMT, given by $[0\ 3.5\ 0\ 0\ 0\ 0]^T$, to $\mathcal{N}_f$ which is generated with an initial condition given by $\bar{X}_0^f = [0\ 0.5\ 0\ 0\ 0\ 0.0051]^T$. Trajectories generated using both a weighted and unweighted virtual net calculated using Theorem 2 are shown. Note again that the trajectory calculated using the weighted adjacency matrix includes more intermediate NMTs, but requires less control (fuel) to execute.

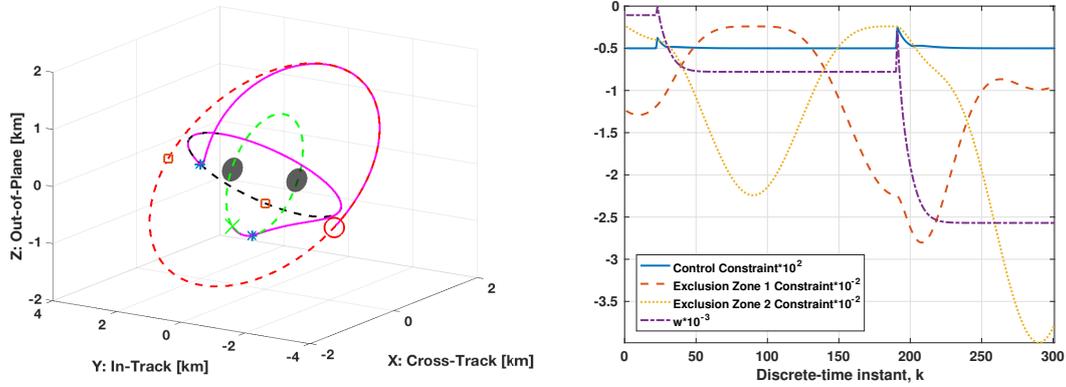
(a) Trajectories from trajectory planned using Procedure 1 and unweighted adjacency matrix

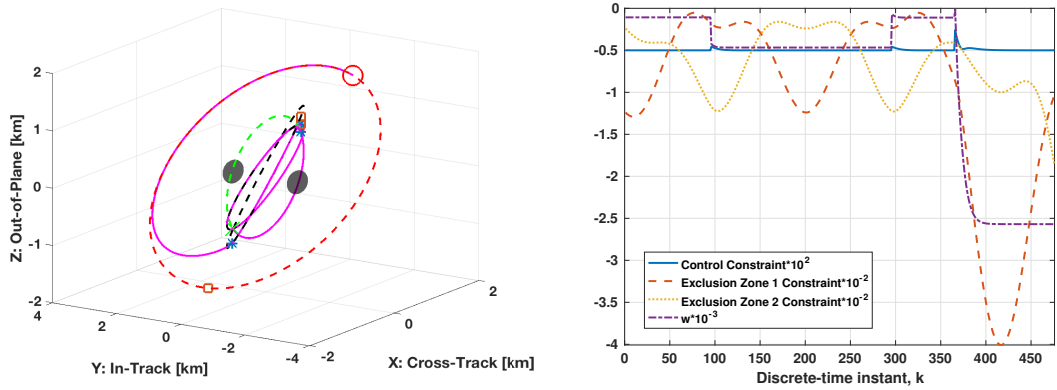
(b) Trajectories from trajectory planned using Procedure 1 and weighted adjacency matrix

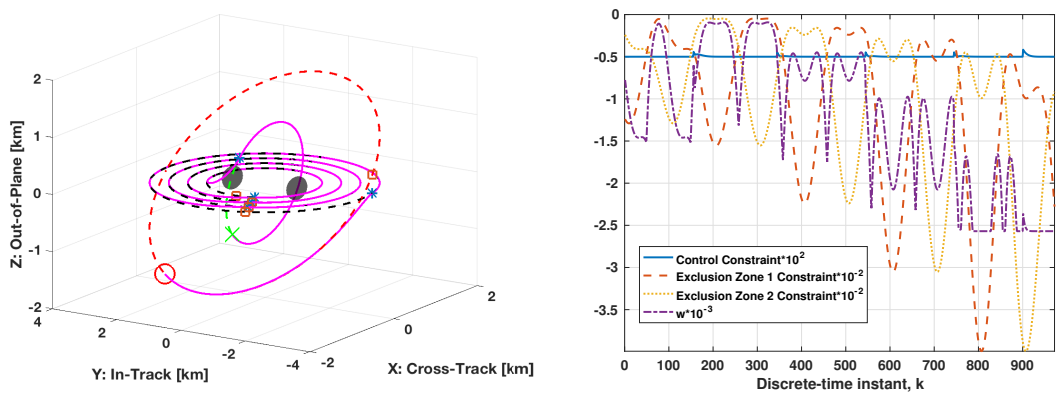
(c) Trajectories from trajectory planned using Procedure 2 and weighted adjacency matrix

Fig. 6   Simulation results using Procedures 1 and 2

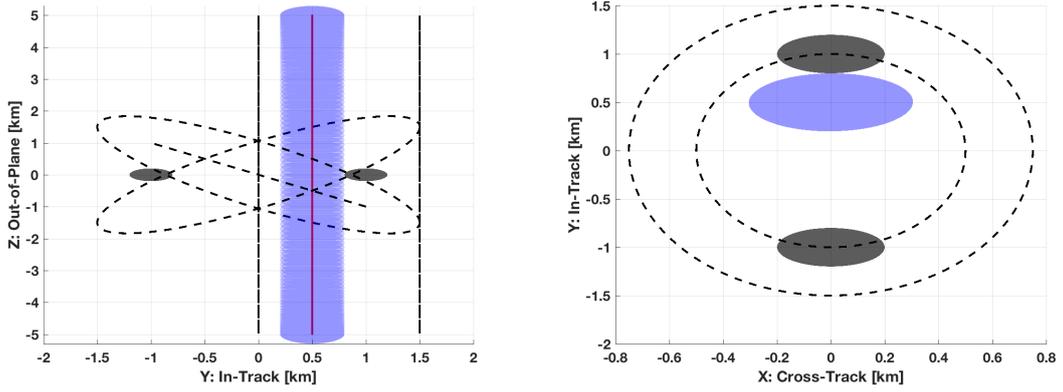

(a) **Two views showing lack of connections to $\mathcal{N}_f$ using the safe, positively-invariant tube calculated using Procedure 1.**

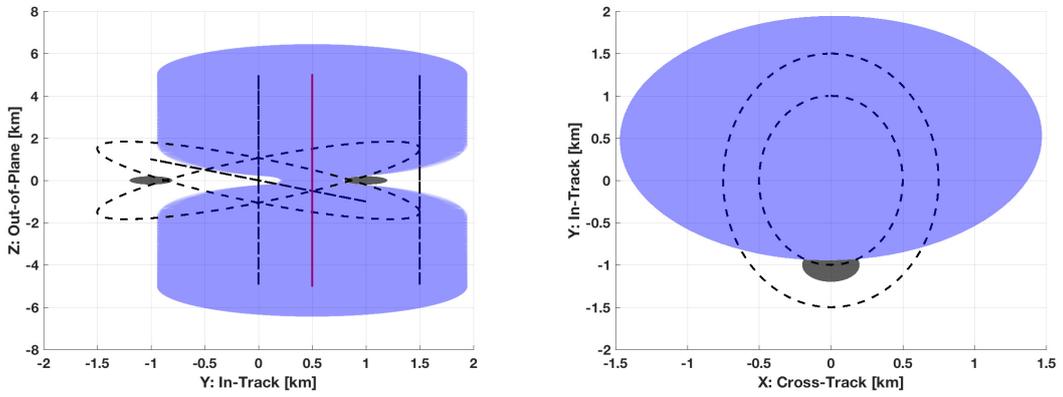

(b) **Two views showing possible connections to $\mathcal{N}_f$ using the positively-invariant tube calculated using Procedure 2.**

**Fig. 7**

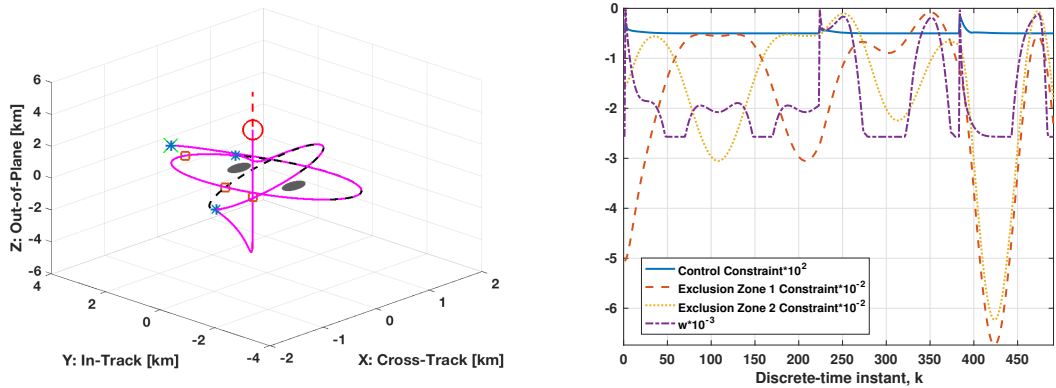

(a) **Trajectories from trajectory planned using unweighted adjacency matrix. Total cost: $J = 2887$ N·sec.**

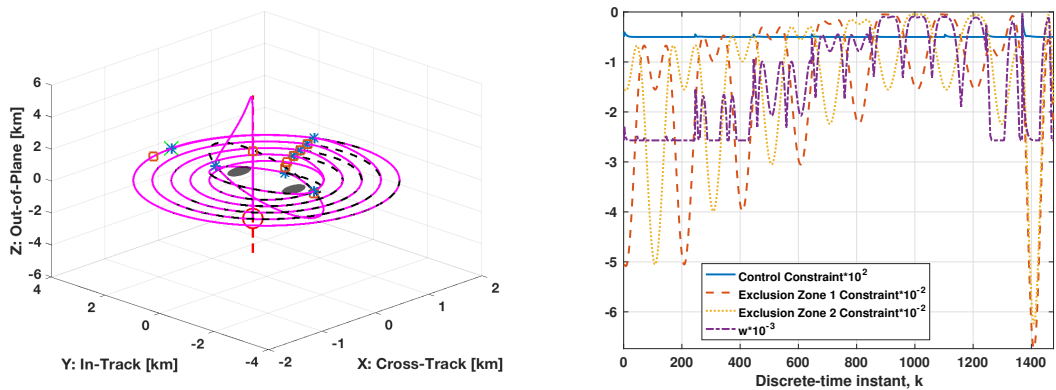

(b) **Trajectories from trajectory planned using weighted adjacency matrix. Total cost: $J = 2290$ N·sec.**

**Fig. 8    Simulation results using Procedure 2**

**Remark 10:** Note that the lack of connections to $\mathcal{N}_f$ using a virtual net calculated using Procedure 1, illustrated in Figure 7a, is partially due to the spacing between NMTs in the virtual net. If more NMTs were added to the virtual net near $\mathcal{N}_f$, then connections may be possible. However, the addition of more NMTs to the virtual net comes at the expense of additional computations required to form the adjacency matrix and connection array.

**Remark 11:** Using the methods described above, for the virtual net considered here with 84 closed NMTs, the approximate computation times required to generate the safe, positively-invariant tube scale factors and the weighted and unweighted adjacency matrices are shown in Table 4. Calculations are performed running MATLAB® R016a on a MacBook Pro® with a 2.8 GHz processor. Note that computation time for the weighted and unweighted connection arrays is negligible since the calculations simply consist of storing indices determined in the adjacency calculation. After forming the virtual net, the calculation of a safe trajectory is accomplished on the order of 0.1 sec. For implementation, the calculations to form the virtual net may be conducted offline, and, after uploading the adjacency matrix, connection array, and NMT initial conditions to the satellite, individual trajectories may be planned on-board. Methods to speed up the offline calculation of the virtual net are currently being investigated.

Table 4   Approximate computation times

|  | Computation Time [min] |
|---|---|
| Safe, Pos.- Inv. Tube Scale Factors using Proc. 1 | 0.8 |
| Safe, Pos.- Inv. Tube Scale Factors using Proc. 2 | 1.5 |
| Unweighted Adjacency Matrix | 13.7 |
| Weighted Adjacency Matrix | 63.0 |

## 6   Conclusion

In this paper, the problem of constrained spacecraft relative motion planning was reduced to a graph search by forming a "virtual net" with nodes corresponding to closed natural motion trajectories. The adjacency of the nodes in the virtual net was determined by forming safe, positively-invariant tubes, that were defined as the union of safe ellipsoidal sets centered at discrete points along the NMT. Two methods to compute the ellipsoidal set scale factors were described and proven to yield safe, positively-invariant tubes. By appropriately weighting the virtual net adjacency matrix, and utilizing a connection array that provided information used to execute safe transfers between NMTs, fuel efficient trajectories were planned using graph search algorithms. Simulation results showed that the developed methodology can be used to generate feasible maneuver solutions to the difficult, non-convex problem of trajectory planning with obstacle avoidance. Similarly to what has been shown in [22] for a virtual net of forced equilibria, it is expected that an additional benefit of our framework is the ability to incorporate both bounded disturbances, such as actuation, navigation and modelling errors, and avoid moving obstacles. These developments are left to future work.

## Appendix A: Calculation of Initial Conditions for Closed NMTs

A "stationary point" closed NMT may be generated with an initial condition $\bar{X}_0$ satisfying

$$y(0) = y_0, \\ x(0) = z(0) = \dot{x}(0) = \dot{y}(0) = \dot{z}(0) = 0, \tag{51}$$

and a "periodic line segment" closed NMT may be generated with an initial condition $\bar{X}_0$ satisfying

$$y(0) = y_0, \quad z(0) = c\sin(\psi), \quad \dot{z}(0) = \omega c \cos(\psi), \\ x(0) = \dot{x}(0) = \dot{y}(0) = 0, \tag{52}$$

where $c \in \mathbb{R}$ gives the magnitude of oscillation, i.e., one-half the length of the line segment, $y_0$ gives the location of intersection with the y-axis, and the phase angle $\psi$ can be arbitrarily chosen.

While elliptical trajectories may be generated centered at any point along the y-axis, trajectories centered at the origin are of particular interest since the origin of Hill's frame is frequently a point of special significance. For example, the origin may be the location of another spacecraft or the center point of a spacecraft formation. A closed elliptical NMT centered at the origin can be generated with any initial condition satisfying (6) and

$$y(0) = \tfrac{2}{\omega}\dot{x}(0). \tag{53}$$

These trajectories can be characterized by three parameters, a scale factor $b$ and two angles, $\theta_1$ and $\theta_2$ [38]. The angles $\theta_1$ and $\theta_2$ are measured from the origin with respect to the relative orbit normal vector, $\hat{h}$, perpendicular to the relative orbital plane, as shown in Figure 9.

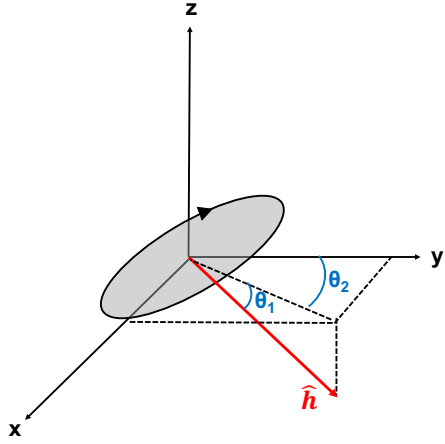

**Fig. 9** Depiction of the angles $\theta_1$ and $\theta_2$ used to parametrize elliptical NMTs centered at the origin

Given $b$, $\theta_1$ and $\theta_2$, the initial condition $\bar{X}_0$ is given by

$$\bar{X}_0 = \begin{bmatrix} b\sin(\nu) & 2b\cos(\nu) & c\sin(\psi) & b\omega\cos(\nu) & -2b\omega\sin(\nu) & c\omega\cos(\psi) \end{bmatrix}^T, \tag{54}$$

where $\nu$ is the $x-y$ plane phase angle corresponding to the initial condition,

$$\tan(\nu - \psi) = 2\tfrac{\cos(\theta_1)}{\tan(\theta_2)}, \tag{55}$$

and

$$c = \tfrac{b}{\sin(\theta_1)}\sqrt{(\tan^2(\theta_2) + 4\cos^2(\theta_1))}. \tag{56}$$

Therefore, by choosing $\nu \in [0, 2\pi]$, and specifying $b$, $\theta_1$ and $\theta_2$, the desired initial condition may be calculated using (54)-(56). Derivations for (54)-(56) can be found in [38].

## Appendix B.1: Proof of Lemma 2

Suppose (38), (39) hold.
Define
$$\begin{aligned} D_1 &= \{e(k+1) \mid e(k+1)^T P e(k+1) = \rho_{k+1}\}, \\ D_2 &= \{e(k) \mid e(k) = \bar{A}^{-1} e(k+1) \mid e(k+1) \in D_1\}, \\ D_3 &= \{e(k) \mid e(k)^T P e(k) = \rho_k\}, \\ D_4 &= \{e(k+1) \mid e(k+1) = \bar{A} e(k) \mid e(k) \in D_3\}, \end{aligned} \tag{57}$$

and consider

$$e'(k+1) \in D_1, \ e'(k) \in D_2, \ e''(k) \in D_3, \text{ and } e''(k+1) \in D_4. \tag{58}$$

Then, by (38), (39),

$$e'(k)^T P e'(k) \geq e'(k+1)^T P e'(k+1) + d(\rho_{k+1}) = \rho_{k+1} + d(\rho_{k+1}) \geq \rho_k. \tag{59}$$

From (59) it follows that any $e''(k)$ can be expressed as $e''(k) = \lambda e'(k)$ where $\lambda \in [0,1]$ if $\rho_{k+1} + d(\rho_{k+1}) = \rho_k$ and $\lambda \in [0,\nu]$, $\nu \in [0,1)$ if $\rho_{k+1} + d(\rho_{k+1}) > \rho_k$. Therefore,

$$e''(k+1)^T P e''(k+1) = e''(k)^T \bar{A}^T P \bar{A} e''(k) = \lambda^2 e'(k)^T \bar{A}^T P \bar{A} e'(k) = \lambda^2 e'(k+1)^T P e'(k+1) \leq \rho_{k+1}, \tag{60}$$

where the last inequality holds for $\lambda \in [0,1]$ and becomes a strict inequality if $\lambda \in [0,\nu]$, $\nu < 1$.

Considering (60),

$$e''(k)^T P e''(k) - e''(k+1)^T P e''(k+1) = \rho_k - e''(k+1)^T P e(k+1) \geq \rho_k - (\rho_{k+1} - \epsilon), \tag{61}$$

where $\epsilon = 0$ if $\lambda \in [0,1]$ and $0 < \epsilon \leq \rho_{k+1}$ if $\lambda \in [0,\nu]$. Hence, $\bar{d}(\rho_k)$ defined in (41) is given by

$$\bar{d}(\rho_k) = \rho_k - (\rho_{k+1} - \epsilon) \tag{62}$$

and

$$\rho_{k+1} = \rho_k + \epsilon - \bar{d}(\rho_k) \geq \rho_k - \bar{d}(\rho_k), \tag{63}$$

therefore (40) holds.

The proof in the opposite direction is similar. Suppose (40) and (41) hold. Define

$$\begin{aligned} E_1 &= \{e(k) \mid e(k)^T P e(k) = \rho_k\}, \\ E_2 &= \{e(k+1) \mid e(k+1) = \bar{A} e(k) \mid e(k) \in E_1\}, \\ E_3 &= \{e(k+1) \mid e(k+1)^T P e(k+1) = \rho_{k+1}\}, \\ E_4 &= \{e(k) \mid e(k) = \bar{A}^{-1} e(k+1) \mid e(k+1) \in E_3\}, \end{aligned} \tag{64}$$

and consider

$$\tilde{e}(k) \in E_1, \ \tilde{e}(k+1) \in E_2, \ \bar{e}(k+1) \in E_3, \ \bar{e}(k) \in E_4. \tag{65}$$

Then, by (40), (41),

$$\tilde{e}(k+1)^T P \tilde{e}(k+1) \leq \tilde{e}(k)^T P \tilde{e}(k) - \bar{d}(\rho_k) = \rho_k - \bar{d}(\rho_k) \leq \rho_{k+1}. \tag{66}$$

From (66) it follows that any $\bar{e}(k+1)$ can be written as $\bar{e}(k+1) = \lambda \tilde{e}(k+1)$, where $\lambda \in [1,\sigma]$, $\sigma \in \mathbb{R}_{>1}$ if $\rho_k - \bar{d}(\rho_k) = \rho_{k+1}$ and $\lambda \in [1+\nu, \sigma]$, $\nu \in (1,\sigma]$ if $\rho_k - \bar{d}(\rho_k) < \rho_{k+1}$. Therefore,

$$\bar{e}(k)^T P \bar{e}(k) = \bar{e}(k+1)^T \bar{A}^{-1T} P \bar{A}^{-1} \bar{e}(k+1) = \lambda^2 \tilde{e}(k+1)^T \bar{A}^{-1T} P \bar{A}^{-1} \tilde{e}(k+1) = \lambda^2 \tilde{e}(k)^T P \tilde{e}(k) \geq \rho_k, \tag{67}$$

where the last inequality holds for $\lambda \in [1,\sigma]$ and becomes a strict inequality if $\lambda \in [1+\nu, \sigma]$. Considering (67),

$$\bar{e}(k)^T P \bar{e}(k) - \bar{e}(k+1)^T P \bar{e}(k+1) = \bar{e}(k)^T P \bar{e}(k) - \rho_{k+1} \geq (\rho_k + \epsilon) - \rho_{k+1}, \tag{68}$$

where $\epsilon = 0$ if $\lambda \in [1,\sigma]$ and $0 < \epsilon \leq \rho_k$ if $\lambda \in [1+\nu, \sigma]$. Hence, $d(\rho_{k+1})$ defined in (39) is given by

$$d(\rho_{k+1}) = (\rho_k + \epsilon) - \rho_{k+1}, \tag{69}$$

and

$$\rho_k = \rho_{k+1} + d(\rho_{k+1}) - \epsilon \leq \rho_{k+1} + d(\rho_{k+1}), \tag{70}$$

hence 38 holds. □

## Appendix B.2: Proof of Theorem 3

Choosing $\rho_k$ using Procedure 1 ensures that the resulting set $\rho_k$, $k \in [0, k_{max}]$ satisfies Theorem 1. Hence, $\mathcal{T}_\mathcal{N}$ is positively-invariant. Because Procedure 1 results in $\rho_k \leq \rho_k^s$ for all $k \in [0, k_{max}]$, $\mathcal{T}_\mathcal{N}$ is also safe. Because $u_{max} > 0$ and $\mathcal{N} \cap \mathcal{O}_i(s_i, S_i) = \emptyset$ for $i = 1, 2, ..., l$, $\min_{k \in [0, k_{max}]} \rho_k^s > 0$. Therefore, using the Procedure 1 to set $\rho_k$ yields $\rho_k > 0$ for all $k \in [0, k_{max}]$. □

## Appendix B.3: Proof of Theorem 4

Choosing $\rho_k$ using Procedure 2 ensures that the resulting set $\rho_k$, $k \in [0, k_{max}]$ satisfies Theorem 2. Hence, $\mathcal{T}_\mathcal{N}$ is positively-invariant. Because Procedure 2 results in $\rho_k \leq \rho_k^s$ for all $k \in [0, k_{max}]$, $\mathcal{T}_\mathcal{N}$ is also safe. It remains to show that using Procedure 2, $\rho_k > 0$ for all $k \in [0, k_{max}]$. Because $u_{max} > 0$ and $\mathcal{N} \cap \mathcal{O}_i(s_i, S_i) = \emptyset$ for $i = 1, 2, ..., l$., $\min_{k \in [0, k_{max}]} \rho_k^s > 0$. Per step 1 of Procedure 2, the first $\rho_k$ is chosen to be $\rho_{k+1} = \min_{k \in [0, k_{max}]} \rho_k^s > 0$. Therefore, to show that the Procedure 2 results in $\rho_k > 0$ for all $k \in [0, k_{max}]$, it suffices to show that $d(\rho_{k+1}) > 0$ for $\rho_{k+1} > 0$. From Assumptions (A1) and (A2),

$$e(k)^T P e(k) - e(k+1)^T P e(k+1) = e(k+1)\bar{A}^{-1} Q \bar{A}^{-1} e(k+1) \geq \lambda_{min}(Q) \|\bar{A}^{-1} e(k+1)\|_2^2 > 0 \quad (71)$$

for $e(k+1) \neq 0$. Therefore, $d(\rho_{k+1}) \geq \lambda_{min}(Q) \|\bar{A}^{-1} e(k+1)\|_2^2 > 0$ for $\rho_{k+1} > 0$. □

## Appendix C: Converting QCQP to LP

The solution to the QCQP

$$d(\rho_{k+1}) = \min_{e(k+1)} e(k)^T P e(k) - e(k+1)^T P e(k+1) \quad \text{subject to } e(k+1)^T P e(k+1) = \rho_{k+1}, \quad (72)$$

is obtained by reformulating the QCQP as a Linear Program (LP) [45], which can be solved efficiently by many direct methods. First, (72) is re-written using the error dynamics,

$$d(\rho_{k+1}) = \min_{e(k+1)} e(k+1)^T \bar{Q} e(k+1) \quad \text{subject to } e(k+1)^T P e(k+1) = \rho_{k+1}, \quad (73)$$

where $\bar{Q} = \bar{A}^{-1T} P \bar{A}^{-1} - P$. The matrices $\bar{Q} > 0$ and $P > 0$ are simultaneously diagonalized with an invertible matrix $L$ such that $L\bar{Q}L^T = I_{6\times 6}$, and $LPL^T = P_D$, where $I_{6\times 6}$ is the $6 \times 6$ identity matrix and $P_D$ is a diagonal matrix [46]. The matrix $L$ is calculated as follows: $L = (TU)^{-1}$, where $T = VD$, the matrix $V$ has columns corresponding to the normalized eigenvectors of the matrix $\bar{Q}$, $D$ is a diagonal matrix with entries consisting of the square roots of eigenvalues of the matrix $\bar{Q}$, and the matrix $U$ has columns corresponding to the normalized eigenvectors of the matrix $T^{-1} P T^{-1}$.

Let $e(k+1) = L^T y$. With this substitution, problem (72) becomes

$$d(\rho_{k+1}) = \min_y y^T y \quad \text{subject to } y^T P_D y - \rho_{k+1} = 0. \quad (74)$$

Next, define $z_i = y_i^2$, where $y_i$ denotes the $i^{th}$ entry of the vector $y$, and let $P_{Di}$ denote the $(i, i)$ entry of $P_D$. Then, the problem (74) is re-stated as an LP,

$$d(\rho_{k+1}) = \min_{z_i, \ i=1,2,...,6} \sum_{i=1}^6 z_i \quad \text{subject to } \sum_{i=1}^6 [P_{Di} z_i] - \rho_{k+1} = 0, \quad z_i \geq 0, \ i = 1, 2, ...6. \quad (75)$$